\documentclass[12pt,fleqn,leqno]{article}
\usepackage{graphicx}
\usepackage{longtable}
\newcommand{\be}{\begin{equation}}
\newcommand{\ee}{\end{equation}}
\newcommand{\ba}{\begin{array}}
\newcommand{\ea}{\end{array}}

\newcommand{\s}{\sigma}

\newcommand{\bt}{\beta}
\newcommand{\de}{\delta}

\renewcommand{\t}{\theta}

\newcommand{\bea}{\begin{eqnarray}}
\newcommand{\eea}{\end{eqnarray}}

\newcommand{\ga}{\gamma}
\newcommand{\w}{\omega}
\newcommand{\ep}{\varepsilon}
\newcommand{\Ga}{\Gamma}
\newcommand{\la}{\lambda}

\newcommand{\pri}{\prime}

\newcommand{\ze}{{\zeta}}
\newcommand{\hy}{{\;_2F_1}}

\begin{document}
\newtheorem{pro}[thm]{Proposition}
\newtheorem{lem}[thm]{Lemma}
\newtheorem{cor}[thm]{Corollary}
\title{Small eigenvalues of large Hankel matrices}
\author{Yang Chen and Nigel Lawrence\\
Department of Mathematics\\
Imperial College\\
180 Queen's Gate\\
London, SW7 2BZ}
%\date{March 17, 1996}

\maketitle

\begin{abstract}
In this paper we investigate the smallest 
eigenvalue, denoted as $\la_N,$  of a $(N+1)\times (N+1)$ Hankel or 
moments matrix, associated  with the weight, $w(x)=\exp(-x^{\bt}),
\;x>0,\;\;\bt>0$, in the large $N$ limit.
Using a previous result, the asymptotics
for the polynomials, $P_n(z),\;z\notin[0,\infty)$, orthonormal with 
respect to $w,$ which are required in the determination of $\la_N$ are found. 
Adopting an argument of Szeg\"{o} the asymptotic behaviour of 
$\la_N$, for $\bt>1/2$ where the related moment problem is determinate,
is derived.
This generalises the result given by Szeg\"{o} for
$\bt=1$. It is shown that for $\bt>1/2$ the smallest eigenvalue
of the infinite Hankel matrix is zero, while for $0<\bt<1/2$ it is greater
then a positive constant. This shows a phase transition in the corresponding
Hermitian random matrix model as the parameter $\bt$ varies with $\bt=1/2$
identified as the critical point. The smallest eigenvalue at 
this point is conjectured.      

\end{abstract}

\bigskip

\setcounter{section}{1}
\setcounter{equation}{0}
\setcounter{thm}{0}

{\bf 1. Introduction.}

In the theory of Hermitian random matrices, the Hankel determinant plays
an important role,
\bea
D_N=\det_{0\leq i,j\leq N}(\mu_{i+j})\;.\nonumber
\eea

For a given weight function $w(t)$ on $J$ ($\subseteq R$,) the
moments $\mu_k$ are 
\bea
\mu_k:=\int_Jw(t)t^{k}dt\;\; ;\;\;k=0,1,2,\cdots
\eea
Associated with $w(t)$ is a Hankel matrix or moment matrix of order $N+1$,
$\{H_{jk}\},$ whose entries are given by
\bea
H_{jk}:=\mu_{j+k}\;\;;\;\;0\leq j,k\leq N\;.
\eea
It is believed that correlations between eigenvalues of random matrices
are universal after a suitable rescaling. In the following treatment we
will show that a fundamental quantity, namely the least eigenvalues 
of these Hankel matrices exhibit a critical dependence on the weight
function. It is this non-universal property that motivates our investigation 
of this problem.     

If $J$ is a single interval say $[a,b]$, where $a$ and $b$ are fixed
and the Szeg\"o condition,
\bea
\int_{a}^{b}\frac{v(x)dx}{{\sqrt {(b-x)(x-a)}}}<\infty,\;\;v:=-\ln w,\nonumber
\eea
is satisfied then the asymptotic behaviour of the Hankel determinants for
large $N$ was established by Szeg\"{o}, \cite{Sz3}.
Let $\la_N$ denote the smallest 
eigenvalue. Szeg\"o also investigated the behaviour of $\la_N$ for 
large $N$ \cite{Sz}. He 
studied the cases for which $J$ can either be a finite or infinite interval
with special choices for $w.$ If $w(x)=1,\;x\in(-1,1)$ and $w(x)=1,\;x\in(0,1),
$ then the respective smallest eigenvalues are for large $N$ \footnote{ 
Throughout this paper, the relation, $a_N\simeq b_N$ means $\lim_{N\to\infty}
a_N/b_N=1.$}
\bea
\mbox{} \la_N&\simeq&2^{\frac{9}{4}}\pi^{\frac{3}{2}}N^{\frac{1}{2}}
({\sqrt 2}-1)^{2N+3}
\nonumber\\
\mbox{} \la_N&\simeq&2^{\frac{15}{4}}\pi^{\frac{3}{2}}N^{\frac{1}{2}}
({\sqrt 2}-1)^{4N+4}.
\nonumber
\eea
Widom and Wilf \cite{ww} generalised Szeg\"o's results to a kind
of ``universal'' law. Thus if $w(x)>0,\;x\in[a,b]$ and 
the Szeg\"o condition is satisfied then it was found in \cite{ww} that
\bea
\la_N\simeq A\;N^{\frac{1}{2}}\;B^{-N},\nonumber
\eea
where $A$ and $B$ are computable constants depending on $w$, $a$, $b$ and 
are independent of $N.$

In \cite{Sz}, Szeg\"o also considered the cases of infinite intervals where 
$w(x)={\rm exp}\left[-x^2\right],\;x\in(-\infty,+\infty)$ and 
$w(x)={\rm exp}\left[-x\right],\;x\in[0,+\infty),$ are the weights
of the Hermite and Laguerre polynomials\footnote{There is a factor of 4
missing from the original formula for $\la_N$; the last equation in page 677 
of \cite{Sz}.}. The respective smallest eigenvalues are
\bea
\mbox{} \la_N&\simeq& 2^{\frac{13}{4}}\pi^{\frac{3}{2}}\;e\;N^{\frac{1}{4}}
{\rm exp}
\left[-2(2N)^{\frac{1}{2}}\right],\nonumber\\  
\mbox{} \la_N&\simeq& 2^{\frac{7}{2}}\pi^{\frac{3}{2}}\;e\;N^{\frac{1}{4}}
{\rm exp}
\left[-4N^{\frac{1}{2}}\right].\nonumber
\eea
Observe that in the examples given above the smallest eigenvalues are 
exponentially small. Therefore it is very hard to numerically invert the 
Hankel matrices associated with these weights.

It is well known that  $\la_N$ is given by the Rayleigh
quotient:
\bea
\la_N=\min\left\{\frac{\sum_{j,k=0}^N H_{jk}x_j\overline{x}_k}{\sum_{j=0}^N
|x_j|^2}\right\}\;.
\eea
If $\pi_N(z)$ is a polynomial of degree $N$, with coefficients
$x_j,\;j=0,...,N$
\bea
\pi_N(z):=\sum_{j=0}^{N}x_jz^j\;,
\eea
then
\bea
\sum_{j,k=0}^N H_{jk}x_j\overline{x}_k=\int_J|\pi_N(t)|^2w(t)dt
\eea
and
\bea
\sum_{j=0}^N|x_j|^2=\int_0^{2\pi}|\pi_N(e^{i\phi})|^2\frac{d\phi}{2\pi}.
\eea
Consequently we can rephrase the extremal expression for $\la_N$, (1.3), as,
\bea
\frac{2\pi}{\la_N}=\max\left\{\int_0^{2\pi}|\pi_N(e^{i\phi})|^2d\phi\;:\;
\int_J|\pi_N(t)|^2w(t)dt=1\right\}\;.
\eea
Let $\{P_n(t)\}$ be the polynomials, orthonormal with respect to  $w(t)$,
then $\pi_N$ has the expansion, 
\bea
\pi_N(z)=\sum_{j=0}^{N}c_jP_j(z)\;.
\eea
Thus 
\bea
\int_0^{2\pi}|\pi_N(e^{i\phi})|^2d\phi=\sum_{j,k=0}^N
K_{jk}c_j\overline{c}_k\;,
\eea
where
\bea
K_{jk}:=\int_0^{2\pi}P_j(z)\overline{P_k(z)}d\phi \;;\;z=e^{i\phi}\;.
\eea
Therefore (1.7) is equivalent to
\bea
\frac{2\pi}{\la_N}=\max\left\{\sum_{j,k=0}^NK_{jk}c_j\overline{c}_k\;:\;
\sum_{j=0}^{N}|c_j|^2=1\right\}\;.
\eea
With the Schwarz inequality, which states that for all values of
$j$ and $k$
\bea
|K_{jk}|\leq K_{jj}^{\frac{1}{2}}K_{kk}^{\frac{1}{2}}\;,\nonumber
\eea
and Cauchy's inequality we obtain an upper bound of (1.11):
\bea
\mbox{}\sum_{j,k=0}^NK_{jk}c_{j}\overline{c}_{k}&\leq&
\sum_{j,k=0}^N|K_{jk}||c_j||c_k|\\
\mbox{}&\leq&\sum_{j,k=0}^N K_{jj}^{\frac{1}{2}}K_{kk}^{\frac{1}{2}}
|c_{j}||c_{k}|
\nonumber\\
\mbox{}&\leq&\left(\sum_{j=0}^N  K_{jj}\right)
\left(\sum_{j=0}^N |c_{j}|^2\right)\nonumber\\
\mbox{}&=&\sum_{j=0}^N K_{jj}.\nonumber
\eea
Therefore a lower bound for the smallest eigenvalue $\la_N$ is given by
\bea
\frac{2\pi}{\sum_{j=0}^{N}K_{jj}}\leq \lambda_N.
\eea

This paper is organised as follows: In section 2, by adopting a previous
result \cite{Ch:La}, we obtain the asymptotic formula for the polynomials
orthonormal with respect to $w(t):=\exp[-t^{\bt}],\;\bt>1/2,$ which is
then employed in sections 3 and 4 for the determination of the large $N$
behaviour of $\la_N.$ In these sections we show, following \cite{Sz},
by an appropriate choice of the vector $\{c_j\}$, that the lower bound given by
(1.13) is in fact an asymptotic estimate for large $N$. By a simple 
application of Laplace's method, $\sum_{j=0}^{N}K_{jj}$ is estimated.
Thus the asymptotic form of $\la_N$ follows. In order to test the
accuracy of the theory, these results are checked against numerical 
calculations for various $\bt$ and $N$, which were obtained using the Jacobi 
rotation algorithm \cite{NR} to reduce the Hankel matrix to diagonal form.
This is found in section 5.

\setcounter{section}{2}
\setcounter{equation}{0}
\setcounter{thm}{0}

{\bf 2. The weight } $w(t)=\exp[-t^{\bt}],\;t\in[0,\infty)$.

In this case, the moments are
\bea
\mu_n=\frac{1}{\bt}\Gamma\left(\frac{n+1}{\beta}\right).
\eea
In order to find a lower bound for the smallest eigenvalue good knowledge 
is required of the associated orthonormal polynomials $\{P_N(z)\}$, 
for $N$ large and $z\notin(0,\infty).$ In \cite{Ch:La}, by applying
the linear statistics formula for matrix ensembles together with the 
Heine's determinant representation, asymptotic forms for the
polynomials with weight
$w(t)={\rm exp}[-v(t)]$, where $v(t)$ is an arbitrary convex function
supported on $[0,\infty)$, are derived. The zeros of these 
polynomials are supported on $(a,b)\subset {\bf R}$. Here $a=0$, 
whilst $b(N)$ follows from the condition that ensures that $P_N(t)$ has $N$ 
roots on $(a,b)$, one finds that \cite{Ch:La},
\bea
b(N;\bt)=CN^{\frac{1}{\bt}},\;\;{\rm where\;}
C=C(\bt):=4\left[\frac{\Gamma^2(\bt)}{\Gamma(2\bt)}\right]^{\frac{1}{\bt}}
N^{\frac{1}{\bt}}\;.
\eea
The normalised 
polynomials as $N\to\infty$ are found, using \cite{Ch:La}, to be
\bea
P_N(t)\simeq\frac{(-1)^N}{\sqrt{2\pi b}}\frac{
\exp[-f(t)+(2N+1)\ln(\sqrt{\ze}+\sqrt{1+\ze})]}{[\ze(1+\ze)]^{\frac{1}{4}}}
,\;\zeta:=-\frac{t}{b},\;\;t\notin[0,b]\;,
\eea
where $f$ is given by
\bea
f(t):=\frac{\sqrt{t(t-b)}}{2\pi}\int_0^b\frac{dy}{y-t}\frac{y^{\bt}}{\sqrt{
y(b-y)}}\;,t\notin[0,b].
\eea
From the definition and basic properties of the hypergeometric functions
\cite{Gr:Ry},
\bea
\mbox{} f(t)&=&-\frac{N}{\bt-\frac{1}{2}}{\sqrt {\ze(1+\ze)}}
\hy\left(1,1-\bt;\frac{3}{2}-\bt;-\ze\right)
- \frac{(-t)^{\bt}}{2}\sec\pi\bt\\
\mbox{} &=&-\frac{N}{\bt}{\sqrt {\frac{\ze}{1+\ze}}}
\hy\left(1,\frac{1}{2};\bt+1;\frac{1}{1+\ze}\right)\;.\nonumber
\eea
At this point note the dichotomy of the problem,  the nature of the
Hypergeometric function dictates that whilst the first representation is
more convenient in the large $b$ limit, where $|\ze|<<1$, it cannot be used
when $\bt=n+\frac{1}{2},\;\;n=1,2,\dots$, necessitating
the use of the second result of (2.5) in such instances.

Using the fact that
\bea
\ln(\sqrt{\ze}+\sqrt{1+\ze})=\sqrt{\ze}\hy\left(\frac{1}{2},\frac{1}{2};
\frac{3}{2};-\ze\right)\;,
\eea
we find,
\bea
(2N+1)\ln(\sqrt{\ze}+\sqrt{1+\ze})\simeq\frac{(-t)^{\bt}}{\sqrt{\pi}C^{\bt}}
\sum_{r=0}^{E[\bt-\frac{1}{2}]}(-1)^ra_r\ze^{r+\frac{1}{2}-\bt}\;,
\eea
where $E[n]$ denotes the integer part of $n$ and
\bea
a_r:=\frac{\Gamma(r+\frac{1}{2})}{(r+\frac{1}{2})\Gamma(r+1)}\;.
\eea
So the asymptotic expression of the
polynomials for $t\notin(0,\infty)$, is,
\bea
P_N(t)\simeq\frac{(-1)^N\ze^{\frac{1}{4}}}{\sqrt{-2\pi t}}\exp
\left(-f(t)+\frac{(-t)^{\bt}}{C^{\bt}\sqrt{\pi}}\sum_{r=0}^{E[\bt-\frac{1}{2}
]}(-1)^ra_r\ze^{r+\frac{1}{2}-\bt}\right)\;.
\eea
To make further progress we now consider separately the 
two possible cases, as identified above, for $\bt>1/2$.

\setcounter{section}{3}
\setcounter{equation}{0}
\setcounter{thm}{0}

{\bf 3.} $\bt\neq n+\frac{1}{2},\;\;n=1,2,3\cdots$

When $\bt\neq n+\frac{1}{2}$, we use the first form for $f(t)$
in equation (2.5). The series expansion for the function $\hy(
{\scriptstyle 1,1-\bt;\frac{3}{2}-\bt;-\ze})$, valid for $|\ze|<1$, is
\bea
\hy\left(1,1-\bt;\frac{3}{2}-\bt;-\ze\right)=
\frac{\Ga(\frac{3}{2}-\bt)}{\Ga(1-\bt)}\sum_{r=0}^{\infty}(-1)^r
\frac{\Ga(1-\bt+r)}{\Ga(\frac{3}{2}-\bt+r)}\ze^r\;,
\eea
whilst for $|\ze|<1$, $\sqrt{1+\ze}$ may be written as
\bea
\sqrt{1+\ze}=\frac{-1}{2\sqrt{\pi}}\sum_{r=0}^{\infty}(-1)^r
\frac{\Ga(r-\frac{1}{2})}{\Ga(r+1)}\ze^r\;.
\eea
With this noted, the expansion for $f(t)$  as $\ze\to 0$ is
\bea
\mbox{}f(t)&\simeq&-\frac{1}{2\sqrt{\pi}}\frac{\Gamma(\frac{1}{2}-\bt)}
{\Gamma(1-\bt)}\left(\frac{-t}{C}\right)^{\bt}
\sum_{r=0}^{E[\bt-\frac{1}{2}]}(-1)^rb_r\ze^{r+\frac{1}{2}-\bt}\\
\mbox{}&-&\frac{(-t)^{\bt}}{2}\sec\pi\bt\;,\nonumber
\eea
where
\bea
b_r:=\sum_{s=0}^r\frac{\Gamma(s-\frac{1}{2})\Gamma(1-\bt+r-s)}
{\Gamma(s+1)\Gamma(\frac{3}{2}-\bt+r-s)}\;.
\eea
Recall that $\ze=-tC^{-1}N^{-\frac{1}{\bt}}$, and by the use of equation 
(2.9) we have,
\bea
\mbox{}P_N(t)&\simeq&\frac{(-1)^N}{\sqrt{2\pi}}
(-tCN^{\frac{1}{\bt}})^{-\frac{1}{4}}\exp\left[\frac{(-t)^{\bt}}{2}
\sec\pi\bt\right]\\
\mbox{}&\times&\exp\left[\frac{N^{1-\frac{1}{2\bt}}}{\sqrt{\pi C}}\sum_{r=0}
^{E[\bt-\frac{1}{2}]}(-1)^rA_r\frac{(-t)^{r+\frac{1}{2}}}
{(CN^{\frac{1}{\bt}})^r}\right]\nonumber\;,
\eea
with
\bea
A_r:=a_r+\frac{\Gamma(\frac{1}{2}-\bt)}{2\Gamma(1-\bt)}b_r\;.
\eea

Note with $\bt=1$, we find $C=4$ and
$A_0=4\sqrt\pi$ and consequently recover the classical result for the
Laguerre polynomials due to Perron \cite{Sz2},
\bea
P_N(t)\simeq\frac{(-1)^N}{2\sqrt{\pi}}(-tN)^{-\frac{1}{4}}
\exp\left[2\sqrt{-tN}+\frac{t}{2}\right]\;\;,\;t\notin[0,\infty),\;
N\to\infty\;.
\eea

With $P_N(t)$ having the form (3.5), where
$A_0=\frac{4\sqrt{\pi}\bt}{2\bt-1}$ is positive for $\bt>1/2$, we observe
that for sufficiently large $j$ and $k$ the dominant contributions to
$K_{jk}$ are from the arc of the
unit circle around $t=-1$. Thus by fixing an arbitrary positive number $\w$
and confining ourselves to values of $j$ and $k$ satisfying
\bea
N-\w N^{\frac{1}{2\bt}}\leq j,k \leq N\;,
\eea
we have
\bea
K_{jk}\simeq\int_{\pi-\ep}^{\pi+\ep}P_j
\left(e^{i\phi}\right)P_k\left(e^{-i\phi}\right)d\phi\;.
\eea
Using the substitution $\t=\phi-\pi$ and expanding the integrand
for $|\t|<<1$ gives the following,
\bea
\mbox{}K_{jk}&\simeq&\frac{(-1)^{j+k}}{2\pi\sqrt{C}}e^{\sec\pi\bt}
N^{-\frac{1}{2\bt}}\\
\mbox{}&\times&
\int_{-\ep}^{\ep}\exp\Biggl[
\frac{1}{\sqrt{\pi C}}\sum_{r=0}^{E[\bt-\frac{1}{2}]}(-1)^r\frac{A_r}{C^r}
\Biggl[\left(1-\frac{(2r+1)^2\t^2}{8}\right)
\left(j^{1-\frac{1}{2\bt}-\frac{r}{\bt}}+
k^{1-\frac{1}{2\bt}-\frac{r}{\bt}}\right)\nonumber\\
\mbox{}&&\hspace{0.5in}
 +\frac{(2r+1)i\t}{2}\left(j^{1-\frac{1}{2\bt}
-\frac{r}{\bt}}-
k^{1-\frac{1}{2\bt}-\frac{r}{\bt}}\right)\Biggr]
\Biggr]d\theta\;.\nonumber
\eea
Because $j^{1-\frac{1}{2\bt}+\frac{r}{\bt}}
 -k^{1-\frac{1}{2\bt}-\frac{r}{\bt}}$
remains bounded in the range specified by (3.8) we can disregard the
linear term in $\t$ in the integrand.
This integral can then be approximated by extending the range of
integration to the real axis, which does not affect the asymptotic behaviour,
as contributions from $(-\infty,-\ep)$ and $(\ep,\infty)$ are sub-dominant
compared to those from $[-\ep,\ep]$ as $j,k\to\infty$. Therefore,
\bea
\mbox{}K_{jk}&\simeq&
\frac{(-1)^{j+k}}{(\pi C)^{\frac{1}{4}}}A_0^{-\frac{1}{2}}e^{\sec\pi\bt}
N^{-\frac{1}{2}-\frac{1}{4\bt}}\\
\mbox{}&\times&\exp\left[\frac{1}{\sqrt{\pi
C}}\sum_{r=0}^{E[\bt-\frac{1}{2}]}
(-1)^r\frac{A_r}{C^r}\left(j^{1-\frac{1}{2\bt}-\frac{r}{\bt}}+
k^{1-\frac{1}{2\bt}-\frac{r}{\bt}}\right)\right]\;.
\nonumber
\eea

From (3.11), we see that when $j$ and $k$ are sufficiently large and satisfy
(3.8), 
\bea
K_{jk}\simeq(-1)^{j+k}K_{jj}^{\frac{1}{2}}K_{kk}^{\frac{1}{2}}\;.
\eea
This is especially useful as it enables the
determination of the large $N$ behaviour of $\la_N$.
By choosing the vector $\{c_j\}$, as in \cite{Sz}, such that
\bea
c_{j}=\cases{\s e^{i\pi j}K_{jj}^{\frac{1}{2}}\;\;&if
$E[N-\w N^{\frac{1}{2\bt}}]\leq j\leq N$ \cr
0&if $j<E[N-\w N^{\frac{1}{2\bt}}]\;,$\cr}
\eea
where $\s$ is a positive number determined by the condition
\bea
\sum_{j=0}^N|c_j|^2=\s^2\sum_{j=E[N-\w N^{\frac{1}{2\bt}}]}^N K_{jj}=1\;,
\eea
we find, using (3.12) and (3.14), that
\bea
\mbox{}\sum_{j,k=0}^N K_{jk}c_j\overline{c}_k&=&
\sum_{j,k=E[N-\w N^{\frac{1}{2\bt}}]}^N
\s^2e^{i\pi(j-k)}K_{jk}K_{jj}^{\frac{1}{2}}K_{kk}^{\frac{1}{2}}\\
\mbox{}&\simeq&\s^2\left[\sum_{j=E[N-\w N^{\frac{1}{2\bt}}]}^N
K_{jj}\right]^2\nonumber\\
\mbox{}&=&\sum_{j=E[N-\w N^{\frac{1}{2\bt}}]}^N K_{jj}\nonumber\;.
\eea
Recalling equation (1.11),
we see that since $\w$ is arbitrarily large the
asymptotic behaviour of the maximum, by virtue of the
inequality (1.13), is well approximated by $\sum_{j=0}^N K_{jj}$. Therefore
we have shown that
\bea
\frac{2\pi}{\la_N}\simeq\sum_{j=0}^NK_{jj}\;.
\eea
The leading behaviour of this sum for large $N$ is in turn found by 
replacing the sum by an integral and by applying Laplace's method, which 
in this context may be stated as :

If for $x\in[a,b]$, the real continuous function $\phi(x)$ has as its
maximum
the value $\phi(b)$, then as $N\to\infty$
\bea
\int_a^b f(x)\exp[N\phi(x)]dx\simeq\frac{f(b)\exp[N\phi(b)]}
{N\phi^{\prime}(b)}\;.
\eea

A simple calculation gives the expression for $\la_N$,
\bea
\mbox{}\frac{2\pi}{\la_N}&\simeq&\frac{1}{4}\pi^{-\frac{1}{4}}C^{\frac{1}{4}
}
A_0^{-\frac{1}{2}}e^{sec\pi\bt}N^{-\frac{1}{2}+\frac{1}{4\bt}}\\
\mbox{}&\times&\exp\left[\frac{2N^{1-\frac{1}{2\bt}}}{\sqrt{\pi C}}
\sum_{r=0}^{E[\bt-\frac{1}{2}]}(-1)^r\frac{A_r}{C^r}N^{-\frac{r}{\bt}}
\right]\;.\nonumber
\eea
Putting $\bt=1$, Szeg\"{o}'s classical result 
for the Laguerre weight is recovered:
\bea
\frac{2\pi}{\la_N}\simeq 2^{-\frac{5}{2}}\pi^{-\frac{1}{2}}e^{-1}
N^{-\frac{1}{4}}\exp[4\sqrt{N}]\;.
\eea
From (3.18) we see that the smallest eigenvalue is exponentially small 
for large $N$ and is zero for the corresponding infinite Hankel matrix.

\setcounter{section}{4}
\setcounter{equation}{0}
\setcounter{thm}{0}

{\bf 4. $\bt=n+\frac{1}{2}\;,\;n=1,2,\cdots $}

In this section we investigate the case where $\bt=n+\frac{1}{2}\;,
\;n\geq 1$.
Such cases, as was explained previously, require the second form of $f(t)$
in (2.5). To obtain the asymptotic expansion for $f(t)$,
we first note the following result for the hypergeometric function :

If $\bt=n+\frac{1}{2}$ with $n=1,2,\dots$ then
\bea
\hy\left(1,\frac{1}{2};\bt+1;x\right)=L_{\bt}\frac{(x-1)^{\bt-\frac{1}{2}}}
{x^{\bt+\frac{1}{2}}}\left(\sqrt{x}\ln\left[\frac{1+\sqrt{x}}{1-\sqrt{x}}
\right]
+\sum_{r=1}^{\bt-\frac{1}{2}}\frac{1}{L_{r-\frac{1}{2}}}
\left(\frac{x}{x-1}\right)^r\right)\;,
\eea
where $L_{r}$ is given by
\bea
L_{r}:=\frac{r}{2\pi}C^{r}(r)\;.
\eea
This is easily be proved by using an inductive argument,
noting the following version of Gauss' recursion relations \cite{Gr:Ry}
\bea
\mbox{}\hy\left(1,\frac{1}{2};n+\frac{5}{2};z\right)&=&
\frac{(n+\frac{3}{2})(z-1)}{(n+1)z}\left[
\hy\left(1,\frac{1}{2};n+\frac{3}{2};z\right)-
\hy\left(1,\frac{1}{2};n+\frac{1}{2};z\right)\right]\\
\mbox{}&+&\frac{n(n+\frac{3}{2})}{(n+1)(n+\frac{1}{2})}\hy\left(
1,\frac{1}{2};n+\frac{3}{2};z\right)\;\nonumber
\eea
together with the fact that
\bea
\hy\left(1,\frac{1}{2};\frac{5}{2};z\right)=\frac{3}{4}\frac{(z-1)}
{z^{\frac{3}{2}}}\ln\left[\frac{1+\sqrt{z}}{1-\sqrt{z}}\right]
+\frac{3}{2}z\;.
\eea
Therefore,
\bea
f(t)=\frac{(-1)^{\bt+\frac{1}{2}}}{2\pi}(-t)^{\bt}\left(
\ln\left[\frac{\sqrt{1+\ze}+1}{\sqrt{1+\ze}-1}\right]
+\sqrt{1+\ze}\sum_{r=1}^{\bt-\frac{1}{2}}(-1)^r
\frac{\ze^{-r}}{L_{r-\frac{1}{2}}}
\right)\;.
\eea
Using (3.2), we find
\bea
\mbox{}f(t)&\simeq&\frac{(-1)^{\bt+\frac{1}{2}}(-t)^{\bt}}{2\pi}
\ln\left[\frac{4}{\ze}\right]\\
\mbox{}&+&\frac{(-t)^{\bt}}{4\pi^{\frac{3}{2}}}\sum_{r=0}^{\bt-\frac{1}{2}}
(-1)^r\de_{\bt-\frac{1}{2}-r}\ze^{r+\frac{1}{2}-\bt}\;,\;|\zeta|<<1\nonumber
\eea
where
\bea
\de_r:=\sum_{s=1}^{\bt-\frac{1}{2}}\frac{\ga_{s-r}}{L_{s-\frac{1}{2}}}
\eea
and
\bea
\ga_r:=\cases{\frac{\Ga(r-1/2)}{\Ga(r+1)}\;\;&if $r\geq0$\cr
0&if $r<0$.\cr}
\eea
Recalling $\ze =-tC^{-1}N^{-\frac{1}{\bt}}$, the
strong asymptotics of the polynomials for $t\notin[0,\infty)$ reads,
\bea
\mbox{}P_N(t)&\simeq&\frac{(-1)^N}{\sqrt{2\pi}}
(-tCN^{\frac{1}{\bt}})^{\frac{1}{4}}\exp\left[\frac{(-1)^{\bt-\frac{1}{2}}
(-t)^{\bt}}{2\pi}\ln\left(\frac{4CN^{\frac{1}{\bt}}}{-t}\right)\right]\\
\mbox{}&\times&\exp\left[
\frac{N^{1-\frac{1}{2\bt}}}{\sqrt{\pi C}}\sum_{r=0}^{\bt-\frac{1}{2}}
(-1)^rB_r\frac{(-t)^{r+\frac{1}{2}}}{(CN^{\frac{1}{\bt}})^r}
\right]\;,\nonumber
\eea
where
\bea
B_r:=a_r-\frac{L_{\bt}}{2\bt}\de_{\bt-\frac{1}{2}-r}\;.
\eea
Note the appearance of the logarithm in exponential.
Since  $B_0=\frac{4\sqrt{\pi}\bt}{2\bt-1}>0$ and using an argument similar
to that in the previous section, we see that in determining $K_{jk}$
the essential contribution comes from the arc in the vicinity of $t=-1$.
As before restricting $j,k$ to the range given in (3.8), we have,
\bea
K_{jk}\simeq\int_{-\ep}^{\ep}P_j(-e^{i\t})P_k(-e^{-i\t})d\t\;.
\eea
We expand the exponential
in the integrand for $|\t|<<1$, keeping terms up to second order
and then extend the range of integration to the infinite interval. Because
$j^{1-\frac{1}{2\bt}-\frac{r}{\bt}}-k^{1-\frac{1}{2\bt}-\frac{r}{\bt}}$ and
$\ln(j/k)$ remain bounded in the range given by (3.8), we find
\bea
\mbox{}K_{jk}&\simeq&\frac{(-1)^{j+k}}{(\pi C)^{\frac{1}{4}}}
B_0^{-\frac{1}{2}}N^{-\frac{1}{2}-\frac{1}{4\bt}}(4CN^{\frac{1}{\bt}})^
{\frac{(-1)^{\bt-\frac{1}{2}}}{\pi}}\\
\mbox{}&\times&\exp\left[
\frac{1}{\sqrt{\pi C}}\sum_{r=0}^{\bt-\frac{1}{2}}(-1)^r\frac{B_r}{C^r}
\left(j^{1-\frac{1}{2\bt}-\frac{r}{\bt}}+k^{1-\frac{1}{2\bt}-\frac{r}{\bt}}
\right)\right]\;.\nonumber
\eea
Again note that for sufficiently large $j$ and $k$, satisfying (3.8),
\bea
K_{jj}\simeq(-1)^{j+k}K_{jj}^{\frac{1}{2}}K_{kk}^{\frac{1}{2}}\;.
\eea
Repeating the argument of the previous section, it follows that
\bea
\frac{2\pi}{\la_N}\simeq\int_0^N K_{jj}dj\;.
\eea
The leading term in the asymptotic expansion of this integral as
$N\to\infty$ follows from an
application of Laplace's method and is given by
\bea
\mbox{}\frac{2\pi}{\la_N}&\simeq&\frac{1}{4}\pi^{-\frac{1}{4}}
C^{\frac{1}{4}}B_0^{-\frac{1}{2}}
N^{-\frac{1}{2}+\frac{1}{4\bt}}(4CN^{\frac{1}{\bt}})
^{\frac{(-1)^{\bt-\frac{1}{2}}}{\pi}}\\
\mbox{}&\times&\exp\left[\frac{2N^{1-\frac{1}{2\bt}}}{\sqrt{\pi C}}
\sum_{r=0}^{\bt-\frac{1}{2}}(-1)^r\frac{B_r}{C^r}N^{-\frac{r}{\bt}}
\right]\;.\nonumber
\eea
Effectively $\exp[\sec\pi\bt]$ in (3.18) is replaced by 
$(4CN^{1/\bt})^{\frac{(-1)^{\bt-1/2}}{\pi}}$. Note the alternating nature
of this additional factor depending on whether $\bt-1/2$ is odd or even.
Again (4.15) shows that $\lim_{N\to\infty}\la_N=0.$ According to standard 
theory\cite{Ak}, the moment problem associated
with $w(x),\;x\geq 0$ is indeterminate if
\bea
\int_{0}^{\infty}\frac{v(x)}{\sqrt{x}(1+x)}dx<\infty.\nonumber
\eea  
Therefore $\bt=1/2$ is special as it marks the transition point at which the
moment problem becomes indeterminate.
Assuming, the result given in (2.9) holds, we have
\bea
P_N(t)\simeq\frac{(-1)^N}{2\pi}(-t)^{-\frac{1}{4}}N^{-\frac{1}{2}}\exp\left[
\frac{\sqrt{-t}}{\pi}\left(\ln\left[\frac{4\pi N}{\sqrt{-t}}\right]
+1\right)\right],\;\;t\notin[0,\infty).
\eea
Again if we confine ourselves to the range where
$j$ and $k$ are sufficiently large to enable the use of the
above asymptotic representation, we find that the major contributions
to $K_{jk}$ are from the arc around $t=-1$. But,
due to the behaviour of $P_N(t)$ with increasing $N$, it is quite
clear that $|K_{jk}|$ decreases as $j,k\to\infty$, making
an analysis analogous to that of the previous sections impossible.

It is however possible to obtain an approximate lower bound
for the least eigenvalue, since (1.13)
still holds. Applying the Christoffel-Darboux formula \cite{Sz2} and
the result given in \cite{ch} for the large $N$ off diagonal recurrence 
coeeficients, we find,
\bea
\mbox{}\sum_{j=0}^NK_{jj}&=&\int_{-\pi}^{\pi}\sum_{j=0}^N
P_{j}(-e^{i\t})P_{j}(-e^{-i\t})d\t\\
\mbox{}&\simeq&\pi^2N^2\int_{-\pi}^{\pi}
\frac{P_{N}(-e^{i\t})P_{N+1}(-e^{-i\t})-
P_{N}(-e^{-i\t})P_{N+1}(-e^{i\t})}{e^{i\t}-e^{-i\t}}d\t\nonumber
\eea
Thus using Laplaces method,
$$\int_{a}^{b}dx\;f(x){\rm exp}[N\phi(x)]
\simeq f(c){\rm exp}[N\phi(c)]{\sqrt {\frac{2\pi}{-Ng^{\pri\pri}(c)}}}\;\;,
\;\;{\rm as}\;\; N\to +\infty,$$
where $c\in(a,b)$ is the maximum of $\phi(x)$ for $x\in(a,b)$, gives
\bea
\sum_{j=0}^NK_{jj}\simeq\frac{(4\pi Ne)^{2/\pi}}{4\sqrt{\ln(4\pi Ne)}}.
\eea
So at the point $\bt=1/2$ the smallest eigenvalue appears to decrease
algebraically instead of exponentially. 

\setcounter{section}{4}
\setcounter{equation}{0}
\setcounter{thm}{0}

{\bf 5. Numerical Results}

In this section we check the accuracy of our asymptotic expressions
for the least eigenvalue of the the various Hankel matrices against
numerical results. Due to the fact that the moment matrices in these cases
are very ill conditioned becuase of the vast range in scale of the matrix
elements, the Jacobi rotation algorithm \cite{NR}, proved far more stable
than the more conventional techniques for numerically determining a small
selection of the eigenvalues of large
symmetric matrices such as the Lanczos procedure or Householder's
method \cite{Go}.  This appears to be an unusual phenomenon.
Because of the behaviour of the
matrix elements in these problems it is necessary to implement a 
multiple-precision package that allows floating point arithmetic of
arbitrary precision. The library of sub-routines created by Brent\cite{Br}
was employed to combat the effect of rounding errors in the numerical 
procedures.  

For $0<\bt<1/2$, 
the corresponding moment problem becomes indeterminate
\cite{Ak}, and as a consequence the sum
\bea
\sum_{j=0}^{\infty}|P_j(z)|^2,\nonumber
\eea
converges for every $z$ in every compact subset of the complex plane. 
Therefore
\bea
\sum_{j=0}^{\infty}K_{jj}=\xi>0,\nonumber
\eea
and the smallest eigenvalue for the corresponding infinite Hankel is a 
positive constant bounded below by $2\pi/\xi$. Proof of the extention 
of the above statement to all indeterminate moment problems and other
related topics can be found in \cite{bci}. The situation for $0<\bt<1/2$ is 
in contrast to the results
for  $\bt>1/2$ where (3.18) and (4.15), as confirmed by the numerics,
show that the sum diverges - A fact that is also well-known from the standard
theory when the moment problem is determinate \cite{Ak}. This separation of 
behaviour in the two regions is the phenomenon of phase transition alluded 
to earlier.     

The comparison between the numerical values of $\la_n$ and those obtained 
from the theoretical expressions (3.18) and (4.15) is shown in table 1. 
and figure 1. below.

\begin{figure}[htb]
\centering
\includegraphics[height=3in]{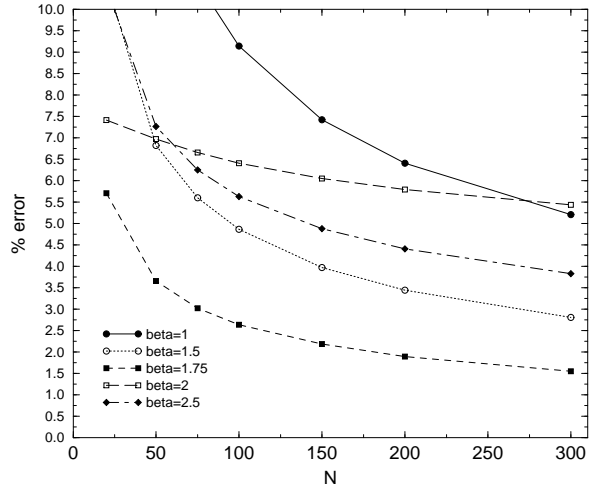}
\caption{The percentage error of the theoretical values of 
$\la_N$ when compared to those obtained numerically, for various $\bt$.}
\end{figure}

\begin{table}
\begin{small}
\begin{center}
\begin{tabular}{|c|c|c|c|} 
\hline
\arraycolsep=1cm

$\bt$ & $N$ & Numerical $\la_N$   & Theoretical $\la_N$     \\ \hline

$1  $ 
      & 50  & $2.0948\times 10^{-10}$& $2.3695 \times 10^{-10}$\\ \hline
      & 100 & $2.1079\times 10^{-15}$& $2.3006 \times 10^{-15}$\\ \hline
      & 150 & $2.9551\times 10^{-19}$& $3.1743 \times 10^{-19}$\\ \hline
      & 200 & $1.6387\times 10^{-22}$& $1.7437 \times 10^{-22}$\\ \hline
      & 300 & $5.5215\times 10^{-28}$& $5.8090 \times 10^{-28}$\\ \hline 
$\frac{3}{2}$ 
      & 50  & $6.4066\times 10^{-22}$& $6.8438 \times 10^{-22}$\\ \hline
      & 100 & $6.2353\times 10^{-36}$& $6.5384 \times 10^{-36}$\\ \hline
      & 150 & $9.9476\times 10^{-48}$& $1.0343 \times 10^{-47}$\\ \hline
      & 200 & $2.8132\times 10^{-58}$& $2.9101 \times 10^{-58}$\\ \hline
      & 300 & $4.6009\times 10^{-77}$& $4.7300 \times 10^{-77}$\\ \hline
$\frac{7}{4}$
      & 50  & $6.4483\times 10^{-27}$& $6.6844 \times 10^{-27}$\\ \hline
      & 100 & $1.6976\times 10^{-45}$& $1.7424 \times 10^{-45}$\\ \hline
      & 150 & $1.5193\times 10^{-61}$& $1.5525 \times 10^{-61}$\\ \hline
      & 200 & $3.9265\times 10^{-76}$& $4.0009 \times 10^{-76}$\\ \hline
      & 300 & $1.4844\times 10^{-102}$& $1.5074\times 10^{-102}$\\ \hline
$2  $ 
      & 50  & $2.7356\times 10^{-31}$& $2.5449 \times 10^{-31}$\\ \hline
      & 100 & $3.8907\times 10^{-54}$& $3.6415 \times 10^{-54}$\\ \hline
      & 150 & $2.9557\times 10^{-74}$& $2.7769 \times 10^{-74}$\\ \hline
      & 200 & $8.9775\times 10^{-93}$& $8.4574 \times 10^{-93}$\\ \hline
      & 300 & $9.5593\times 10^{-127}$&$9.0396 \times 10^{-127}$\\ \hline
$\frac{5}{2}$ 
      & 50  & $2.2384\times 10^{-38}$& $2.4010 \times 10^{-38}$\\ \hline
      & 100 & $1.2580\times 10^{-68}$& $1.3288 \times 10^{-68}$\\ \hline
      & 150 & $5.3195\times 10^{-96}$& $5.5789 \times 10^{-96}$\\ \hline
      & 200 & $1.2155\times 10^{-121}$&$1.2691 \times 10^{-121}$\\ \hline
      & 300 & $1.5236\times 10^{-169}$&$1.5819 \times 10^{-169}$\\ \hline
\end{tabular}
\end{center}
\end{small}
\caption{ Numerical and theoretical values of $\la_N$ for
various $\bt$}
\end{table}

\noindent
e-mail: y.chen@ic.ac.uk, n.lawrence@ic.ac.uk

\end{document}